\documentclass[11pt,a4paper]{article}
\usepackage{amsmath}
\usepackage{amssymb}
\usepackage{latexsym}

\newtheorem{theorem}{Theorem}[section]

\newtheorem{proposition}[theorem]{Proposition}
\newtheorem{definition}[theorem]{Definition}
\newtheorem{corollary}[theorem]{Corollary}

\newtheorem{exmp}[theorem]{Example}
\newtheorem{exmps}[theorem]{Examples}
\newtheorem{rem}[theorem]{Remark}

\newenvironment{examples}{\begin{exmps}\rm}{\end{exmps}}
\newenvironment{remark}{\begin{rem}\rm}{\end{rem}\rm}
\newcommand{\prf}{{\em Proof}. }
\newcommand{\qed}{\hspace*{\fill}$\Box$}

\newcommand{\beeq}[1]{\begin{eqnarray}\label{#1}}
\newcommand{\eneq}{\end{eqnarray}}

\newcommand{\ka}{{\cal A}}

\newcommand{\kc}{{\cal C}}

\newcommand{\kh}{{\cal H}}
\newcommand{\kk}{{\cal K}}

\newcommand{\IC}{{\mathbb C}}

\newcommand{\IP}{{\mathbb P}}

\newcommand{\IR}{{\mathbb R}}
\newcommand{\IZ}{{\mathbb Z}}

\newcommand{\vol}{{\rm vol}}

\newcommand{\verylongarrow}[1]{\hbox to #1{\rightarrowfill}}

\begin{document}

{\parindent0mm{\Large\bf Products of harmonic forms and rational curves}}

\bigskip

\bigskip

{\parindent0mm{\large\bf Daniel Huybrechts\\}
{\small Mathematisches Institut\\
Universit\"at zu K\"oln\\
Weyertal 86-90\\
50931 K\"oln, Germany\\}

\bigskip

\bigskip

\bigskip

\bigskip

The product of closed forms is closed again. The analogous
statement for harmonic forms, however, fails. A priori, there is no reason
why the product of harmonic forms should be harmonic again. This phenomenon
was recently studied by Merkulov \cite{Merk}. He shows that
it leads to a natural $A_\infty$-structure on the cohomology of a K\"ahler
manifold. In the context of mirror symmetry Polishchuk made use of
(a twisted version of) this $A_\infty$-structure on elliptic curves to
confirm Kontsevich's homological version of mirror symmetry in this case
\cite{Poli}.

\bigskip

In this paper we show that this failure of harmonicity in fact happens
quite frequently. It usually is related to certain geometric properties
of the manifolds and to the existence of rational  curves in particular.

Let us briefly indicate the main results for the special case of compact
Ricci-flat K\"ahler manifolds. For a K\"ahler class $\omega\in H^2(X,\IR)$
on such a manifold there exists a unique Ricci-flat K\"ahler form
$\tilde\omega$ representing it. Let $\kh^2(\tilde\omega)$ denote the space
of two-forms harmonic with respect to $\tilde\omega$. Of course, for a
different K\"ahler class $\omega'$ and the representing
Ricci-flat K\"ahler form $\tilde\omega'$ this space might be different.

The main technical result (Prop.\ \ref{ProdKX}) says that
$\kh^2(\tilde\omega)$ is independent of $\omega$ if and only if the top
exterior power of any harmonic $\alpha\in\kh^2(\tilde\omega)$ is again
harmonic. This can be used to interprete the failure of harmonicity of the
top exterior power geometrically. Prop.\ \ref{posconeKX} asserts that there
always exist harmonic two-forms with non-harmonic top exterior power,
whenever the K\"ahler cone (or ample cone) does not form a connected
component of the (integral) cone of all classes $\alpha\in
H^{1,1}(X,\IR)$  with $\int_X\alpha^N>0$.

Note that there are many instances where the K\"ahler cone is strictly
smaller. E.g.\ this is the case for any Calabi-Yau manifold that is
birational to a different Calabi-Yau manifold.

In Sect.\ \ref{K3} we apply the result for K3 surfaces. One finds that on
any K3 surfaces containing a rational curve there exists a harmonic two-form
$\alpha$ such that $\alpha^2$ is not harmonic. This can be extended to
arbitrary K3 surfaces by using the existence of rational curves on nearby
K3 surfaces.

\section{Preparations} 

Let $X$ be a compact K\"ahler manifold. Then $\kk_X\subset H^{1,1}(X,\IR)$
denotes the K\"ahler cone, i.e.\  the open set of all K\"ahler classes
on $X$. For a class $\alpha\in H^{1,1}(X,\IR)$ we usually denote by
$\tilde\alpha\in\ka^{1,1}(X)_\IR$ a closed real $(1,1)$-form representing
$\alpha$. Let us recall the following version of the Aubin-Calabi-Yau
theorem \cite{Besse}

\begin{theorem}\label{CY} ---
Let $X$ be an $N$-dimensional compact K\"ahler manifold
with a given volume form $\vol\in\ka^{N,N}(X)_\IR$. For any K\"ahler class
$\omega\in\kk_X$ there exists a unique K\"ahler form $\tilde
\omega\in\ka^{1,1}(X)_\IR$ representing $\omega$, such that
$\tilde\omega^N=c\cdot\vol$, with $c\in\IR$.
\end{theorem}

Since $\tilde\omega^N$ is harmonic with respect to $\tilde\omega$, this can
be equivalently expressed by saying that any K\"ahler class $\omega$ can
uniquely be represented by a K\"ahler form $\tilde\omega$ with respect to
which the given volume form is harmonic.
Note that the constant $c$ can be computed as $c=\int_X\omega^N/\vol(X)$.

\begin{definition}--- For a given volume form $\vol\in\ka^{N.N}(X)_\IR$ we
let $\tilde\kk_X\subset \ka^{1,1}(X)_\IR$ be the set of K\"ahler forms
$\tilde\omega$ with respect to which $\vol$ is harmonic.
\end{definition}

By the Aubin-Calabi-Yau theorem the natural projection
$\tilde\kk_X\to\kk_X$ is bijective. But, in general, $\tilde\kk_X$ is not an
open subset of a linear subspace of $\ka^{1,1}(X)$ (cf.\ \ref{KXlinear}).
Let $\tilde\omega\in\tilde\kk_X$. The tangent space of $\tilde\kk_X$ at
$\tilde\omega$ can be computed as follows. Firstly, we may write
$\tilde\kk_X=\IR_+\times\tilde\kk_X^c$, where
$\tilde\kk_X^c=\{\tilde\omega\in\kk_X|\tilde\omega^N=c\cdot\vol\}$.
Secondly, the infinitesimal deformations of $\tilde\omega$ in the direction
of $\tilde\kk_X^c$ are of the form $\tilde\omega+\varepsilon\tilde v$,
where $\tilde v$ is a closed real $(1,1)$-form and such that 
$(\tilde\omega+\varepsilon\tilde v)^N=\tilde\omega^N$. The latter
condition gives $\tilde\omega^N+{N\choose
2}\varepsilon\tilde\omega^{N-1}\tilde v=\tilde\omega^N$, i.e.\ $\tilde v$
is primitive.
As any closed primitive $(1,1)$-form is harmonic, this shows that the
tangent space of $\tilde\kk_X^c$ at $\tilde\omega$ is the space
$\kh^{1,1}(\omega)_{\IR,prim}$ of real $\tilde\omega$-primitive
$\tilde\omega$-harmonic $(1,1)$-forms. Thirdly, the $\IR_+$-direction
corresponds to the scaling of $\tilde\omega$ and this tangent direction is
therefore canonically identified with $\IR\tilde\omega$. Altogether, one
obtains that $T_{\tilde\omega}\tilde\kk_X=\kh^{1,1}(\tilde\omega)_\IR$
is the space of real $\tilde\omega$-harmonic $(1,1)$-forms. In particular,
$\tilde\kk_X$ is a smooth connected subset of $\ka^{1,1}(X)_\IR$.
To make this approach rigorous, one completes $\ka^{1,1}(X)$ in the
$L^2$-topology. The projection of the closed forms to cohomology is a 
differential map (use e.g.\ Hodge theory, cf.\ \cite{Dem}). The lifted
K\"ahler cone $\tilde\kk_X$ is the intersection of the space of closed
$L^2$-forms with the space of sections of the submanifold of the bundle of
$(1,1)$-forms that consists of those positive
forms whose top exterior power equals (a scalar multiple of)
the given $(N,N)$-form at every point. 


\begin{definition}---
Let $X$ be a compact K\"ahler manifold with a given volume form. Then one
associates to a given K\"ahler class $\omega\in\kk_X$ the space
$\kh^{p.q}(\omega):=\kh^{p,q}(\tilde \omega)$ of $(p,q)$-forms that are harmonic
with respect to the unique $\tilde\omega\in\tilde\kk_X$ representing $\omega$.
\end{definition}

Note that two different K\"ahler forms $\tilde\omega_1$ and
$\tilde\omega_2$ representing the same K\"ahler class
$\omega_1=\omega_2$ always have different spaces of harmonic $(1,1)$-forms.
Indeed, $\tilde\omega_1$ and $\tilde\omega_2$ are
$\tilde\omega_1$-harmonic respectively
$\tilde\omega_2$-harmonic. Since any class, in particular
$\omega_1=\omega_2$, is represented by a unique harmonic form and
$\tilde\omega_1\ne\tilde\omega_2$, this yields
$\kh^{1,1}(\tilde\omega_1)\ne\kh^{1,1}(\tilde\omega_2)$.
One might ask more generally what the relation is between the spaces of
harmonic forms with respect to different K\"ahler forms not
representing the same K\"ahler class.
It is quite interesting to observe that the dependence of
$\kh^{1,1}(\tilde\omega)$ on the K\"ahler class $\omega$ is related to the
problem discussed in the introduction. We will try to make this more
explicit in the next section.

%
%


\section{How `harmonic' depends on the K\"ahler form}

Let us begin with the following fact which relates the shape of
$\tilde\kk_X$ to the dependence of $\kh^{1,1}(\omega)$ on $\omega$.

\begin{proposition}\label{KXlinear}---
The subspace $\kh^{1,1}(\omega)\subset \ka^{1,1}(X)$ is independent
of $\omega$ if and only if $\tilde\kk_X$ spans an $\IR$-linear
subspace of dimension $h^{1,1}(X)$.
\end{proposition}

\prf Let $\kh^{1,1}(\omega)\subset \ka^{1,1}(X)$ be independent of
$\omega\in\kk_X$. Since for any $\omega\in\kk_X$ the unique
$\tilde\omega\in\tilde\kk_X$ representing it is $\tilde\omega$-harmonic,
 the assumption immediately yields
$\tilde\kk_X\subset \kh^{1,1}(\omega)_\IR$ for any $\omega\in\kk_X$.

Conversely, if $\tilde\kk_X$ spans an $\IR$-linear subspace of dimension 
$h^{1,1}(X)$, then this subspace coincides with the tangent space of
$\tilde\kk_X$ at every point $\tilde\omega\in\tilde\kk_X$. But the latter
was identified with $\kh^{1,1}(\omega)_\IR$. Hence, the linear subspace equals
$\kh^{1,1}(\omega)_\IR$ for any $\omega\in\kk_X$ and $\kh^{1,1}(\omega)$, 
therefore, does not depend on $\omega$.
\qed

\bigskip

\begin{remark}---
The assertion might be rephrased from a slightly different point of view 
as follows. The bijective map $\tilde\kk_X\to\kk_X$ can be used to define
a differentiable map $\kk_X\to\ka^2(X)$ (in the
$L^2$-topology). The proposition then just says
that this map is linear if and only if the Gauss map is constant.
It might be instructive to rephrase
some of the results later on in this spirit, e.g.\ Prop.\ \ref{posconeKX}.
\end{remark}

The next proposition states that the `global' change of $\kh^{1,1}(\omega)$
for $\omega\in\kk_X$ is determined by the `harmonic' behaviour with
respect to a single $\omega\in\kk_X$.

\begin{proposition}\label{ProdKX}---
Let $X$ be a compact K\"ahler manifold of dimension $N$ with a fixed
K\"ahler form $\tilde\omega_0$ and volume form $\tilde\omega_0^N/N!$. Then
the following statements are equivalent:

{\it i)} The linear subspace $\kh^{1,1}(\omega)\subset \ka^{1,1}(X)_\IR$ does
not depend on $\omega\in\kk_X$.

{\it ii)} For all $\alpha\in \kh^{1,1}(\omega_0)$ one has
$\alpha^N\in\kh^{N,N}(\omega_0)$. 
\end{proposition}

\prf Let us assume {\it i)}. By the previous lemma the lifted K\"ahler cone
$\tilde\kk_X$ spans the $\IC$-linear subspace
$\kh^{1,1}(\omega_0)$. Since $\tilde\kk_X$ is open in
$\kh^{1,1}(\omega_0)_\IR$ and all
$\alpha\in\tilde\kk_X$ satisfy the $\IC$-linear equation

\begin{equation}
\alpha^N=(\int_X\alpha^N/\int_X\omega_0^n)\cdot\omega_0^N\label{gl}
\end{equation}

which is an algebraic condition, in fact all $\alpha\in
\kh^{1,1}(\omega_0)$ satisfy (\ref{gl}). Hence, for all
$\alpha\in\kh^{1,1}(\omega_0)$ the top exterior power $\alpha^N$ is
harmonic, i.e.\ {\it ii)} holds true.

Let us now assume {\it ii)}. If $\alpha\in\kh^{1,1}(\omega_0)$, such that
its cohomology class $\omega:=[\alpha]$ is a K\"ahler class, let
$\tilde\omega\in\tilde\kk_X$ denote the distinguished representing K\"ahler
form of $\omega$. If $\alpha$ itself is strictly positive definite, then
the unicity of $\tilde\omega$ and {\it ii)} imply $\alpha=\tilde\omega$.
Thus, the intersection of the closed subset $\kh^{1,1}(\omega_0)_\IR$ with the
open cone of strictly positive definite real $(1,1)$-forms is contained in
$\tilde\kk_X$. This intersection is non-empty, as it contains
$\tilde\omega_0$.
Since $\tilde\kk_X$ is a closed connected subset of this
open cone of the same dimension as $\kh^{1,1}(\omega_0)_\IR$ this yields
$\tilde\kk_X\subset\kh^{1,1}(\omega_0)_\IR$.
By Prop.\ \ref{KXlinear} one concludes that $\kh^{1,1}(\omega)$ does not
depend on $\omega\in\kk_X$.
\qed


\section{The positive cone}

The next proposition is a first step towards a geometric understanding of
the failure of harmonicity of $\alpha^N$ for a harmonic form $\alpha$. To
state it we recall the following notation.

\begin{definition}---
For a compact K\"ahler manifold $X$ the positive cone $\kc_X\subset
H^{1,1}(X,\IR)$ is the connected component  of $\{\alpha\in
H^{1,1}(X,\IR)~|~\int_X\alpha^N>0\}$ that contains the K\"ahler cone.
\end{definition}

Note that by definition $\kk_X\subset\kc_X$. 

\begin{proposition}\label{posconeKX}---
If $X$ is a compact K\"ahler cone such that $\kk_X$ is strictly smaller
than $\kc_X$, then for any K\"ahler form $\tilde\omega$ there exists a
$\tilde\omega$-harmonic $(1,1)$-form $\alpha$ such that $\alpha^N$ is not
$\tilde\omega$-harmonic.
\end{proposition}

\prf Assume that there exists a K\"ahler form $\tilde\omega_0$ such that
for all $\alpha\in\kh^{1,1}(\tilde\omega_0)$ also $\alpha^N$ is
$\tilde\omega_0$-harmonic. We endow $X$ with the volume form
$\tilde\omega_0^N/N!$. 
By Prop.\ \ref{ProdKX} the lifted K\"ahler cone
$\tilde\kk_X$ is contained in $\kh^{1,1}(\tilde\omega_0)$.
Since $\kk_X$ is strictly smaller than $\kc_X$ there exists a sequence
$\omega_t\in\kk_X$ converging towards a $\omega\in\kc_X\setminus\kk_X$.
As $\tilde\kk_X$ is contained in the finite-dimensional space
$\kh^{1,1}(\tilde\omega_0)$ the lifted K\"ahler forms
$\tilde\omega_t\in\tilde\kk_X$ will converge towards a form (!) and not
only a current
$\tilde\omega\in\kh^{1,1}(\tilde\omega_0)\setminus\tilde\kk_X$. As
a limit of strictly positive definite forms $\tilde\omega$ is
still semi-positive definite. Moreover, $\tilde\omega$ is strictly
positive definite at $x\in X$ if and only if
$\tilde\omega^N$ does not vanish at $x$. By assumption
$\tilde\omega^N=c\tilde\omega_0^N$ with
$c=\int_X\omega^N/\int_X\omega_0^N$.
Since $\omega\in\kc_X$, the scalar factor $c$ is strictly positive.
Hence, $\tilde\omega^N$ is everywhere non-trivial. Thus $\tilde\omega$ is
strictly positive definite. This yields the contradiction.\qed

\bigskip

The interesting thing here is that the proposition in particular can be
used to determine the positivity of a class with positive top
exterior power just by studying the space of harmonic forms with
respect to a single given, often very special K\"ahler form:

\begin{corollary}---
Let $X$ be a compact K\"ahler manifold with a given K\"ahler form
$\tilde\omega_0$. If for all $\tilde\omega_0$-harmonic $(1,1)$-forms
$\alpha$ the top exterior power $\alpha^N$ is also $\tilde\omega_0$-harmonic,
then any class $\omega\in\kc_X$ is a K\"ahler class.\qed
\end{corollary}

We conclude this section with a few examples, where the assumption of the
corollary is met {\it a priori}. In the later sections we will discuss
examples where $\kk_X$ is strictly smaller than $\kc_X$ and where Prop.\
\ref{posconeKX} can be used to conclude the `failure' of harmonicity.

\begin{examples} ---
{\it i)} If $X$ is a complex torus and $\omega$ is a flat K\"ahler form,
then harmonic forms are constant forms and their products are again
constant, hence harmonic. In particular, one recovers the fact that on
a torus the K\"ahler cone and the positive cone coincide.

{\it ii)} If  for two K\"ahler manifolds $(X,\tilde\omega)$
and $(X',\tilde\omega')$ with  $b_1(X)\cdot b_1(X')=0$
the top exterior power of any harmonic
$(1,1)$-forms on $X$ or on $X'$ is again harmonic, then the same holds
for the product $(X\times X',\tilde\omega\times\tilde\omega')$.
The additional assumption on the Betti-numbers 
is necessary as the product of two curves shows. Indeed,
any $\varphi\in H^{1,0}(X)$, for a curve $X$, is harmonic, but 
$\varphi\wedge\bar\varphi$ is not. Hence, $\alpha=\varphi\times\bar\varphi
+\bar\varphi\times\varphi$ is a harmonic $(1,1)$-form on $X\times X'$ with
non-harmonic $\alpha^2$.

{\it iii)} If $X$ is a K\"ahler manifold, such that $\kh^{1,1}(\omega)$
does not depend on $\omega$, then the same holds for any smooth
finite quotient of $X$.

{\it iv)} For hermitian symmetric spaces of compact type
it is known that the space of harmonic forms equals the space of
forms invariant under the real form. As the latter space is
invariant under products, the K\"ahler cone of an irreducible
hermitian symmetric space coincides with the positive cone.

\end{examples}


\section{K3 surfaces}\label{K3}

As indicated earlier the behaviour of the K\"ahler cone is  closely
related to the geometry of the manifold. 
We shall study this in more detail for K3 surfaces. The next
proposition follows directly from the description of the K\"ahler cone of a
K3 surface. 

\begin{proposition}---
Let $X$ be a K3 surface containing a smooth rational curve. Then for any
K\"ahler form $\tilde\omega$ there exists an $\tilde\omega$-harmonic form
$(1,1)$-form $\alpha$ such that $\alpha^2$ is not harmonic.
\end{proposition}

\prf If $X$ contains a smooth rational curve, then $\kk_X$ is strictly
smaller than $\kc_X$ and we apply Prop.\ \ref{posconeKX}.
Indeed, a smooth rational
curve $C\subset X$ determines a $(-2)$-class $[C]$, whose perpendicular
hyperplane $[C]^\perp$ cuts $\kc_X$ into two parts and $\kk_X$ is
contained in the part that is positive on $C$.\qed

\bigskip


If the harmonicity of the top exterior powers fails for a K\"ahler manifold
with a given K\"ahler form $(X,\tilde\omega)$ then it should do so for any
small deformation of $(X,\tilde\omega)$. For a Ricci-flat K\"ahler structure
on a K3 surface the argument can be reversed and one can use the existence
of rational curves on arbitrarily near deformations to prove the above
proposition on any K3 surface. 

\begin{corollary}--- Let $X$ be an arbitrary K3 surface.
If $\tilde\omega$ is any hyperk\"ahler form on $X$, then there exists a
$\tilde\omega$-harmonic $(1,1)$-form $\alpha$ such that $\alpha^2$
is not $\tilde\omega$-harmonic.
\end{corollary}

\prf
Let $H^0(X,\Omega^2_X)=\IC\sigma$. Then
$$\begin{array}{rcl}
\kh^2(\tilde\omega)&=&\kh^{1,1}(\tilde\omega)
\oplus\kh^{2,0}(\tilde\omega)\oplus\kh^{0,2}(\tilde\omega)\\
&=&\kh^{1,1}(\tilde\omega)\oplus\IC\sigma\oplus\IC\bar\sigma\\
\end{array}$$
As the space of harmonic forms only depends on the underlying hyperk\"ahler
metric $g$, $\kh^{1,1}(\tilde\omega)\oplus\IC\sigma\oplus\IC\bar\sigma$
contains $\kh^{1,1}(\tilde\omega_{aI+bJ+cK})$ for all
$(a,b,c)\in S^2$. Here, $I,J,$ and $K$ are the three complex structures
associated with the hyperk\"ahler metric $g$.

Assume $\alpha^2$ is $g$-harmonic for all
$\alpha\in\kh^{1,1}(\tilde\omega)$. Since
$\sigma=\tilde\omega_J+i\tilde\omega_K$ (up to a scalar factor) and
since the product of a harmonic form with the K\"ahler form is again
harmonic, also $\sigma\bar\sigma$ is harmonic. This implies that $\alpha^2$
is harmonic for all $\alpha\in\kh^2(\omega)$, as
$\sigma^2=\bar\sigma^2=\alpha\sigma=\alpha\bar\sigma=0$ for
$\alpha\in\kh^{1,1}(\tilde\omega)$. Thus, $\alpha^2$ is $g$-harmonic for all
$\alpha\in\kh^{1,1}(\tilde\omega_{aI+bJ+cK})$ and all $(a,b,c)\in S^2$.
On the other hand, it is well-known that for a non-empty (dense)
subset of $S^2$ the K3 surface $(X,aI+bJ+cK)$ contains a smooth rational
curve. Indeed, if $e\in H^2(X,\IZ)$ is any $(-2)$-class, then the subset 
of the moduli space of marked K3 surfaces for which $e$ is of type $(1,1)$ is
a hyperplane section. This hyperplane section, necessarily, cuts the complete
curve given by the base $\IP^1=S^2$ of the twistor family. Hence, on one of the
K3 surfaces $(X,aI+bJ+cK)$ the class $e$ represents a smooth rational curve.
Contradiction.\qed

\bigskip

\begin{remark}---
What are the bad harmonic $(1,1)$-forms? Certainly $\tilde\omega^2$
is harmonic and for any harmonic form $\alpha$ also $\tilde\omega\alpha$
is harmonic. So, if there is any bad harmonic $(1,1)$-form there must
be also one that is $\tilde\omega$-primitive.
Most likely, it is even true that the square of any primitive harmonic
form is not harmonic. The proof of it should closely follow the arguments
in the proof of Prop.\ \ref{posconeKX}, but there is a slight subtlety
concerning the existence of sufficiently many $(-2)$-classes, that I
cannot overcome for the moment. We sketch the rough idea:
%
%
Assume there exists a $\tilde\omega$-harmonic $\tilde\omega$-primitive real
$(1,1)$-form $\alpha$ such that $\alpha^2$ is $\tilde\omega$-harmonic.
As a $\tilde\omega$-harmonic $\tilde\omega$-primitive $(1,1)$-form
$\alpha$ is also of type $(1,1)$ with respect to any complex structure 
$\lambda=aI+bJ+cK$ induced by the hyperk\"ahler metric corresponding to
$\tilde\omega$ (see Prop.\ 7.5 \cite{Huy}). Moreover, $\alpha$ is also
primitive with respect to all K\"ahler forms $\tilde\omega_\lambda$.
{\it Assume that there exists a complex structure
$\lambda\in S^2$, such that $\kc_X\cap\IR[\alpha]\oplus\IR\omega_\lambda$
is not contained in $\kk_X$.} This condition
can be easily rephrased in terms of $(-2)$-classes and thus becomes a
question on the lattice $3U\oplus2(-E_8)$. It looks rather harmless,
but for the time being I do not know a complete proof of it. 
Under this assumption, we may even assume that in fact $\lambda=I$.
Since $\alpha^2$ is harmonic, in fact $\beta^2$ is harmonic for all
$\beta\in\IR\alpha\oplus\IR\tilde\omega\subset\kh^{1,1}(\omega)$.
Going back to the proof of Prop.\ \ref{posconeKX},
we see that the second part of
it can be adapted to this situation and shows that $\psi^{-1}(\kk_X\cap
\IR[\alpha]\oplus\IR\omega)\subset\IR\alpha\oplus\IR\tilde\omega$,
where $\psi:\tilde\kk_X\to\kk_X$. The space  $\psi^{-1}(\kk_X\cap
\IR[\alpha]\oplus\IR\omega)$ is the space of the distinguished K\"ahler
forms whose classes are linear combinations of $[\alpha]$ and
$\omega$. Therefore, all these forms are harmonic
and linear combinations of $\alpha$ and $\tilde\omega$ themselves.
To conclude, we imitate the proof of Prop.\ \ref{posconeKX} and choose a sequence
$\omega_t\in\kk_X\cap\IR[\alpha]\oplus\IR\omega$ converging towards
$\omega'\in\kc_X\setminus\kk_X$. The corresponding sequence
$\tilde\omega_t\in\tilde\kk_X$ is contained in
$\IR\alpha\oplus\IR\tilde\omega$ and converges towards a form(!)
$\tilde\omega'$.
As in the proof of Prop.\ \ref{posconeKX} this leads to a contradiction.
\end{remark}

\section{Hyperk\"ahler manifolds}

We will try to improve upon Prop.\ \ref{posconeKX} in the case of hyperk\"ahler
manifolds. In particular, we will replace the question whether the top
exterior power $\alpha^N$ of an harmonic form $\alpha$ is harmonic by the
corresponding question for the square of $\alpha$. The motivation
for doing so stems from the general philosophy that hyperk\"ahler manifolds
should be treated in almost complete analogy to K3 surfaces
and that instead of the top intersection pairing one should
consider the Beauville-Bogomolov \cite{Beauv} form as the higher dimensional 
analogue of the intersection pairing for K3 surfaces.

Let us begin by recalling some notations and basic facts. By a compact
hyperk\"ahler manifold $X$ we understand a simply-connected compact
K\"ahler manifold, such that $H^0(X,\Omega^2)=\IC\sigma$, where $\sigma$ is
an everywhere non-degenerate holomorphic two-form. A Ricci-flat K\"ahler
form  $\tilde\omega$ turns out to be a hyperk\"ahler form (cf.\
\cite{Beauv}), i.e.\ there exists a metric $g$ and three complex structures
$I$, $J$, and $K:=IJ$, such that the corresponding K\"ahler forms
$\tilde\omega_{aI+bJ+cK}$ are closed for all $(a,b,c)\in S^2$, such that
$I$ is the complex structure defining $X$, and such that
$\tilde\omega=\tilde\omega_I$. One may renormalize $\sigma$, such that
$\sigma=\tilde\omega_J+i\tilde\omega_K$. In particular, multiplying with
$\sigma$ maps harmonic forms to harmonic forms, for this holds true for
the K\"ahler forms $\tilde\omega_J$ and $\tilde\omega_K$.

The positive cone $\kc_X\subset H^{1,1}(X,\IR)$ is a connected component of
$\{\alpha\in H^{1,1}(X,\IR)~|~q_X(\alpha)>0\}$, where $q_X$ is the
Beauville-Bogomolov form (cf.\ \cite{Beauv,Huy}).

\begin{proposition}---
Let $X$ be a $2n$-dimensional compact hyperk\"ahler manifold with
a fixed hyperk\"ahler form $\omega_0$ and the unique holomorphic two-form
$\sigma$. Then $\alpha^2(\sigma\bar\sigma)^{n-1}$ is harmonic for all
$\alpha\in \kh^{1,1}(\omega_0)$ if and only if the linear subspace
$\kh^{1,1}(\omega)\subset\ka^{1,1}(X)$ does not depend on $\omega\in\kk_X$.
\end{proposition}

\prf Assume that  for all $\alpha\in \kh^{1,1}(\omega_0)$ also
$\alpha^2(\sigma\bar\sigma)^{n-1}$ is harmonic. If $\alpha$ is in
addition strictly positive definite and $\tilde\omega\in\tilde\kk_X$ with
$[\alpha]=\omega$, then
$\alpha^2(\sigma\bar\sigma)^{n-1}=\tilde\omega^2(\sigma\bar\sigma)^{n-1}$.
We adapt Calabi's classical
argument to deduce that in this case $\alpha=\tilde\omega$:
If
$\alpha^2(\sigma\bar\sigma)^{n-1}=\tilde\omega^2(\sigma\bar\sigma)^{n-1}$,
then
$(\alpha-\tilde\omega)(\alpha+\tilde\omega)(\sigma\bar\sigma)^{n-1}=0$.
Since $\alpha$ and $\tilde\omega$ are strictly positive definite, also
$(\alpha+\tilde\omega)$ is  strictly positive definite.
By Lemma 6.1 of \cite{Huy} also
$(\alpha+\tilde\omega)(\sigma\bar\sigma)^{n-1}$
is strictly positive. As $[\alpha]=\omega=[\tilde\omega]$, the difference
$\alpha-\tilde\omega$ can be written as $dd^c\varphi$ for some real function 
$\varphi$. But by the maximum principle the equation
$(\alpha+\tilde\omega)(\sigma\bar\sigma)^{n-1}dd^c\varphi=0$
implies $\varphi\equiv const$. Hence, $\alpha=\tilde\omega$.

As in the proof of Prop.\ \ref{posconeKX} this shows that the intersection of the
closed subset $\kh^{1,1}(\omega_0)_\IR$ with the open cone of strictly
positive definite forms in $\ka^{1,1}(X)_\IR$ is contained in $\tilde\kk_X$
and one concludes that $\tilde\kk_X\subset\kh^{1,1}(\omega_0)_\IR$.

Hence, $\kk_X$ spans a linear subspace of the same dimension and, by Lemma
\ref{KXlinear}
this shows that $\kh^{1,1}(\omega)$ is independent of $\omega\in\kk_X$.

Conversely, let $\kh^{1,1}(\omega)$ be independent of
$\omega\in\kk_X$. Then
$\tilde\kk_X\subset \kh^{1,1}(\omega)_\IR$ for any
$\omega\in\kk_X$. Therefore,
$\alpha^2(\sigma\bar\sigma)^{n-1}=c(\sigma\bar\sigma)^{n}$ with $c\in\IR$
for $\alpha$ in the Zariski-dense open subset
$\tilde\kk_X\subset\kh^{1,1}(\omega)_\IR$. Hence,
$\alpha^2(\sigma\bar\sigma)^{n-1}$ is harmonic for any
$\alpha\in\kh^{1,1}(\omega)$ (cf.\ proof of Prop.\ \ref{ProdKX}). 
\qed

\bigskip


Of course, as for K3 surfaces one expects that $\kh^{1,1}(\omega)$
does in fact depend on $\omega$. This would again follow from the existence
of rational curves in every nearby hyperk\"ahler manifold.
But it would actually be more interesting to reverse the argument:
Assume that $X$ is a hyperk\"ahler manifold, such that for any small
deformation $X'$ of $X$ the K\"ahler cone $\kk_{X'}$ equals $\kc_{X'}$.
I expect that this is equivalent to saying that $\kh^{1,1}(\omega)$
does not depend on $\omega$. If for some other reason than the existence
of rational curves as used in the K3 surface case this can be excluded, then
one could conclude that there always is a nearby deformation $X'$ for
which $\kk_{X'}$ is strictly smaller than $\kc_{X'}$. The latter is
expected to imply the existence of rational curves on $X'$. Along these
lines one could try to attack the Kobayashi conjecture, as the existence of
rational curves on nearby deformations would say that $X$ itself cannot be
hyperbolic. Unfortunately, I cannot carry this through even for K3
surface.


\section{Various other examples}

Here we collect a few examples where algebraic geometry predicts the
failure of harmonicity of the top exterior power of harmonic two-forms. In
all examples this is linked to the existence of rational curves.

\bigskip

{\bf Varieties of general type.} Let $X$ be a non-minimal
smooth variety of general type. As I learned from Keiji Oguiso this
immediately implies that the K\"ahler cone is strictly smaller than the
positive cone. His proof goes as follows:
By definition the canonical divisor $K_X$ is big and by the Kodaira Lemma
(cf.\ \cite{KMM}) it can therefore be written as the sum $K_X=H+E$
of an ample divisor $H$ and an effective divisor $E$ (with rational
coefficients). Consider the segment $H_t:=H+tE$ with $t\in [0,1)$.
If all $H_t$ were contained in the positive cone $\kc_X$, then $K_X$ would
be in the closure of $\kc_X$. If the K\"ahler cone coincided with the
positive cone $\kc_X$, then $K_X$ would be nef, contradicting the
hypothesis that $X$ is not minimal.
Hence $t_0:=\sup\{t|H_t\in\kc_X\}\in(0,1)$. If $H_{t_0}$ is not nef,
then $\kk_X$ is strictly smaller than $\kc_X$. Thus, it suffices to show
that $H_{t_0}$ is not nef. If $H_{t_0}$ were nef then all
expressions of the form $H^{N-i}_{t_0}.H^{i-1}.E$ would be non-negative.
Then
$0=H^N_{t_0}=H^{N-1}_{t_0}(H+t_0E)=H^{N-1}_{t_0}.H+t_0H^{N-1}_{t_0}.E$,
so both summands must vanish. In particular,
$0=H^{N-1}_{t_0}.H=H^2.H^{N-2}_{t_0}+t_0H.H^{N-2}_{t_0}.E$ Again this
yields the vanishing of both terms and in particular $0=H^2.H^{N-2}_{t_0}$.
By induction we eventually obtain $0=H^{N-1}.H_{t_0}$ and, furthermore,
$0=H^{N-1}.H_{t_0}=H^N+t_0H^{N-1}.E$. But this time $H^N>0$ yields the
contradiction. Therefore, for a non-minimal variety of general type
one has $\kk_X\ne\kc_X$ and hence there exist harmonic
(with respect to an arbitrary K\"ahler metric) two-forms with non-harmonic
top exterior power. Note that a non-minimal variety contains rational curves.
As the reader will notice, the above proof goes through on an
arbitrary manifold $X$ that admits a big, but not nef line bundle
$L$ (replacing the canonical divisor). Also in this case the positive
cone and the K\"ahler cone differ.

\bigskip

{\bf Birational Calabi-Yau.}
Let $X$ and $X'$ be birational Calabi-Yau manifolds, i.e.\
$K_X$ and $K_{X'}$ are trivial, then the birational map extends
to an isomorphism or there exist harmonic $(1,1)$-forms on $X$, such that
their top exterior power is not harmonic.
Again, a non-trivial birational correspondence produces rational curves.
As one expects for hyperk\"ahler manifolds that $\kh^{1,1}(\omega)$ does
depend on the hyperk\"ahler form even when $X$ does not contain a rational
curve, e.g.\ for K3 surfaces, it would be interesting to see an example of
a simply-connected Calabi-Yau manifold (in particular not a torus), where
it does not.

The same argument could be applied to the case of different
birational minimal models (minimal models are not unique!).
This shows that in the previous example the K\"ahler cone could be strictly
smaller than the positive cone, even when $K_X$ is nef or ample.

\bigskip

{\bf Blow-ups.} This example is very much in the spirit  of  
the previous two. Let $X$ be a non-trivial blow-up of a  
projective variety $Y$. Then $\kk_X$ is strictly smaller than $\kc_X$ and,
therefore, for any K\"ahler structure on $X$ there exist harmonic
$(1,1)$-forms with non-harmonic maximal exterior power. Indeed, if $L$ is
an ample line bundle on $Y$ then $f^*(L)$ is nef, but not ample, and it is
contained in the positive cone. Hence, $f^*(L)\in\kc_X\setminus\kk_X$.
%
%
Note that also the first example could be proved along these lines.
By evoking the contraction theorem one shows that
any non-minimal projective variety $X$ admits a non-trivial contraction to
a projective variety $Y$. The above argument then yields that $\kk_X$ and
$\kc_X$ are different.

\section{Chern forms}\label{chernforms}

Let $X$ be a compact K\"ahler manifold with a Ricci-flat K\"ahler form
$\tilde\omega$. If $F$ denotes the curvature of the Levi-Cevita connection
$\nabla$, then the Bianchi identity reads $\nabla F=0$. The
K\"ahler-Einstein condition implies $\Lambda_\omega F=0$. The last equation
can be expressed by saying that $F$ is $\tilde\omega$-primitive. Analogously
to the fact that any closed primitive $(1,1)$-form is in fact harmonic, one
has that for $F$ with $\nabla F=0$ the primitivity condition
$\Lambda_{\tilde\omega}F=0$ is equivalent to the harmonicity condition
$\nabla * F=0$. As for untwisted harmonic $(1,1)$-forms one might ask for
the harmonicity of the product $F^m$. Slightly less ambitious, one could
ask whether the trace of this expression, an honnest differential
form, is harmonic. This trace is, in fact,
a scalar multiple of the Chern character $ch_m(X,\tilde\omega)$. 

\bigskip

{\bf Question.} --- Let $(X,\tilde\omega)$ be a Ricci-flat K\"ahler
manifold. Are the Chern forms $ch_m(X,\tilde\omega)$ harmonic with respect
to $\tilde\omega$ ?

\bigskip

By what was said about K3 surface we shall expect a negative answer to this
question at least in this case:

\bigskip

{\bf Problem.} --- Let $X$ be a K3 surface with a hyperk\"ahler form
$\tilde\omega$. Let $c_2\in A^2(X)$ be the associated Chern form. Show that
$c_2$ is not harmonic with respect to $\tilde\omega$ !

\bigskip

So, this should be seen in analogy to the fact that $\alpha^2$ is not
harmonic for any primitive harmonic $(1,1)$-form $\alpha$. Here, $\alpha$
is replaced by the curvatue $F$ and $\alpha^2$ by $tr F^2$. It is likely
that the non-harmonicity of $c_2$ can be shown by standard methods in
differential geometry, in particular by using the fact that $c_2$ is
essentially $\|F\|\cdot \tilde\omega^2$ (see \cite{Besse}), but I do not know
how to do this.

Furthermore, it is not clear to me what the relation between the above
question and  the one treated in the previous sections is.
I could imagine that the non-harmonicity of $ch_m$ in fact implies the
existence of harmonic $(1,1)$-forms with non-harmonic top exterior power.

\bigskip

\bigskip

{\bf Acknowledgements.} I wish to thank U.\ Semmelmann for his interest in
this work and M.\ Lehn and D.\ Kaledin for
making valuable comments on a first version of
it. I am most grateful to Keiji Oguiso for its enthusiastic
help with several arguments.

\bigskip


\bigskip

\bigskip


\end{document}